\input amstex

\def\bhar{B(H^2(m))}
\def\mfi{M\sb{\phi}}
\def\hinfty{H\sp\infty(m)}
\def\falfa{\{F\sb\alpha\}\sb{\alpha\in\Gamma}}
\def\fijalfa{{\phi\sb{j,\alpha}}}
\def\tjalfa{T\sb{j,\alpha}}
\def\jalfa{J\sb{\alpha}}
\def\salfa{\Cal S\sb{\alpha}}
\def\calf{\Cal F}
\def\cstar{$C^*$-algebra\ }

\def\ck{C(K)}
\def\mfialfa{M\sb{\phi\sb{j,\alpha}}}
\def\tilf{\tilde{\Cal F}}
\def\calck{\Cal C(C(K))} 
\def\calinf{\Cal C(L\sp\infty(m))} 
\def\taulinf{\Cal T(L\sp\infty(m))} 
\def\tauck{\Cal T(C(K))} 
\def\har{H^2(m)}
\def\lar{L^2(m)}
\def\linf{L\sp\infty(m)}
\def\tfi{T\sb{\phi}}

\def\caltf{\Cal T(\Cal F)}

\def\calm{\Cal M}
\def\calr{\Cal R}
\def\p00{\Phi\sb{00}}
\def\tj1{T\sb{j\sb 1,\alpha\sb 1}}
\def\tjm{T\sb{j\sb m,\alpha\sb m}}
\def\ti1{T\sb{i\sb 1,\beta\sb 1}}
\def\tin{T\sb{i\sb n,\beta\sb n}}
\def\ttj1{\widehat{T}\sb{j\sb 1,\alpha\sb 1}}
\def\ttjm{\widehat{T}\sb{j\sb m,\alpha\sb m}}
\def\tti1{\widehat{T}\sb{i\sb 1,\beta\sb 1}}
\def\ttin{\widehat{T}\sb{i\sb n,\beta\sb n}}
\def\rhoh{\rho\sb {\Cal H}}
\def\rhok{\rho\sb {\Cal K}}
\def\pik{\pi\sb {\Cal K}}
\def\pk{P\sb {\Cal K}}

\def\jalfa{{J\sb\alpha}}
\def\pix{\pi(X)}
\def\sw{$\sigma$-weakly }
\def\smj{\sum\sb{j\in J\sb\alpha}}
\def\smf{\sum\sb{j\in F}}
\def\ttiljs{\hat{T}\sb{j,\alpha}\sp *}
\def\tjas{T\sb{j,\alpha}\sp *}
\def\tilx{\hat{X}}
\def\tilc{{\Cal C}}
\def\cle{\Cal E}
\def\cstare{C^*(\Cal E)}
\def\bk{{\Cal B}({\Cal K})}

\def\pti{\{T\sb j\}\sb{j\in J}}
\def\calf{{\Cal F}}
\def\cals{{\Cal S}}
\def\salf{\{{\Cal S}\sb\alpha\}\sb{\alpha\in\Gamma}}
\def\salfa{{\Cal S}\sb\alpha}
\def\tja{T\sb{j,\alpha}}
\def\calt{{\Cal T}}

\def\cltf{{\Cal T}({\Cal F})}
\def\tilh{\widehat{\Cal H}}
\def\tilf{\widehat{\Cal F}}
\def\ttilj{\widehat{T}\sb{j,\alpha}}
\def\tsalf{\{\widehat{\Cal S}\sb\alpha\}\sb{\alpha\in 
\Gamma}}
\def\tsalfa{\widehat{\Cal S}\sb\alpha}
\def\btilh{{\Cal B}(\widehat{\Cal H})}
\def\cstartf{C^*({\Cal T}({\Cal F}))}
\def\cstarf{C^*({\Cal F})}
\def\bh{{\Cal B}({\Cal H})}
\def\cstar{$C^*$-algebra }


\def\1{\text{\bf 1}}
\def\B{{\Cal B}}

\def\tr{\text{\rm tr}}
\def\Ker{\text{\rm Ker}\,}
\def\Ran{\text{\rm Ran}\,}

\def\Hc{{\Cal H}}
\def\Kc{{\Cal K}}

\documentstyle{amsppt}
\topmatter
\title Some Exact Sequences
for  Toeplitz Algebras of 
 Spherical 
Isometries 
\endtitle
\author Bebe Prunaru\endauthor
\address Institute of Mathematics ``Simion Stoilow'' of the 
Romanian
Academy,
P.O. Box 1-764, RO-014700 Bucharest, Romania\endaddress
\email  
Bebe.Prunaru\@imar.ro\endemail
\abstract 
A family $\{T_j\}_{j\in J}$ of commuting Hilbert space 
operators is said 
to be a spherical isometry if
$\sum_{j\in J}T^*_jT_j=1$ in the weak operator topology. 
We show that every commuting family 
$\Cal F$ of spherical isometries
has a commuting normal extension
$\hat{\Cal F}.$
Moreover, if $\hat{\Cal F}$ is minimal,
then there exists
a natural  short exact sequence
$0\to\Cal C\to C^*(\Cal F)\to C^*(\hat{\Cal F})\to 0$
with a completely isometric cross-section,
where $\Cal C$ is the commutator ideal in $C^*(\Cal F).$
We also show that the space of Toeplitz operators associated to $\Cal F$ is completely isometric to the commutant of the minimal normal extension $\hat{\Cal F}.$ Applications of these results are given for Toeplitz operators on strictly pseudoconvex or bounded symmetric domains.   
  
\noindent{\it 2000 Mathematics Subject Classification}: 
Primary 47L80; 47B35 
Secondary 47B20; 46L07

\noindent{\it Key words and phrases}: Toeplitz operators; 
Toeplitz algebras; spectral inclusion;   
spherical isometry; completely positive projection; injective operator spaces; 

\endabstract

\leftheadtext{Bebe Prunaru}
\rightheadtext{Exact sequences for Toeplitz algebras}
\endtopmatter

\document

\head 1. Introduction\endhead

A spherical isometry on a Hilbert space $\Hc$  
is a commuting family ${\Cal S}=\{T_j\}_{j\in J}$ 
of bounded operators 
on $\Hc$ such that $\sum\limits_{j\in J}T_j^*T_j=\1$, 
To each spherical isometry one can associate 
its set of  Toeplitz-type  operators consisting of all solutions
 of the operator equation  
$$\sum\limits_{j\in J}T_j^*XT_j=X.$$
One defines in a similar way Toeplitz operators associated to 
arbitrary commuting families of spherical isometries.

In the present paper we apply some 
operator space  techniques in order 
to construct exact sequences  
for  $C^*-$algebras generated by spherical 
 isometries or by their associated  Toeplitz-type 
operators. For this purpose we make use
of  the more or less known fact that the set 
$\calt(\calf)$ of all Toeplitz operators associated to an  
arbitrary family $\calf$ of commuting spherical isometries
is an injective operator system, which means 
that it is the range of a completely positive unital
mapping $\Phi:\bh\to\bh$ with $\Phi\circ\Phi=\Phi.$
This fact, when combined with a certain basic property of such 
projections proved in [CE77],
enables us to show  (see Theorem 2.9) that $\calf$
 admits a commuting normal extension
$\tilf$ on some Hilbert space $\tilh$ containing $\Hc$  
and that if $\tilf$ is minimal then there exists a *-representation $\pi$  
from the   $C^*$-algebra generated by $\calt(\calf)$ 
onto the commutant of $\tilf$ in $B(\tilh)$ for which the 
compression map $\rho(X)=P_{\Hc}X|_{\Hc}$ is a complete isometric cross section whose
 range coincides with $\calt(\calf).$
When we restrict $\pi$ to the \cstar 
$\cstarf$ generated by $\calf$ in $\bh$ we obtain 
a *-representation of $\cstarf$ onto $C^*(\tilf)$
whose kernel coincides with the commutator ideal of $\cstarf.$
We also show that any operator in the commutant of $\calf$ has
 a unique norm-preserving extension to an operator in the commutant of $\tilf.$

It turns out that several classes of Toeplitz operators on 
various Hardy spaces can be realized as common fixed points 
of some commuting families of completely positive mappings 
induced by spherical isometries. In this sense, apart from the 
well-known Brown-Halmos  characterization of Toeplitz 
operators on the  
unit circle  (see \cite{BH63}) there is a similar result due to 
A.M. Davie and N. 
Jewell \cite{DJ77} for the case of Toeplitz operators on the 
unit sphere in ${\Bbb C}^n$ where the unilateral shift is 
replaced by  
the Szego $n-$tuple. Similar characterizations 
 also hold for  Toeplitz operators 
on Hardy spaces  on ordered groups
\cite{Mu87}. 
 The method we shall develop here 
allows us to enlarge considerably the class of Toeplitz 
operators that admit such characterization. 
We show that if $\Omega\subset {\Bbb C}^n$ is either a 
bounded strictly pseudoconvex domain or a bounded symmetric domain     
and $m$ is  any Borel probability measure 
on the Shilov boundary of $\Omega$ then the Toeplitz operators on the corresponding Hardy space $H^2(m)$ are indeed the fixed points of a certain spherical isometry. As a consequence we obtain exact sequences and spectral inclusion theorems for  operators of the type mentioned above.

\head 2. Spherical isometries and their associated Toeplitz operators 
\endhead

We recall for later use the following by-product of the proof of Theorem~3.1 
in \cite{CE77}, which is also stated as Lemma~6.1.2 in \cite{ER00}. 

\proclaim{Theorem~2.1}
Let $\Phi \colon\bh\to\bh$ be a completely positive and 
completely contractive mapping such that $\Phi ^2=\Phi $. 
Then for all $X,Y\in \bh$ we have
$$\Phi(\Phi(X)Y)= \Phi(X\Phi(Y))= \Phi(\Phi(X)\Phi(Y)).$$
\endproclaim

In \cite{CE77} and \cite{ER00} 
this result is used to show that the range of $\Phi$ is 
completely  isometric to a \cstar, where the 
multiplication is defined by the rule
$$\Phi (X)\circ \Phi (Y)=\Phi (\Phi (X) \Phi (Y))$$
for every $X,Y\in\bh$. 
We shall recover the latter 
result as a consequence of Theorem~2.2 below.
To begin with, let us fix some notation.
Let $\Phi \colon \bh\to\bh$ be a  completely positive  and 
completely contractive mapping such that $\Phi ^2=\Phi $ 
and let $\cle=\Ran \Phi $ denote its range. 
We denote by 
$\cstare$ the unital \cstar generated by $\cle$ in 
$\bh$.  
Let
$$\Phi _0\colon \cstare\to\bh$$ denote the restriction of $\Phi $ 
to $\cstare$.
 Then, according to the Stinespring dilation theorem 
(see for instance Theorem~5.2.1 in \cite{ER00}), 
there exist a Hilbert space ${\Kc}$, a bounded operator  
$V\colon {\Hc}\to {\Kc}$ and a unital *-representation  
$$\pi\colon \cstare\to\bk$$ such that 
$\Phi _0(X)=V^*\pi(X)V$ for all $X\in\cstare$.
Thus the diagram 
$$\CD
C^*({\Cal E}) @>{\pi}>> \bk \\
@V{\Phi _0}VV @ VV{\rho}V \\
{\Cal E}  @>{\iota}>> \bh
\endCD
$$
is commutative, where $\rho(X)=V^*XV$ whenever 
$X\in\bk$, 
and $\iota$ is the inclusion map.  
We shall assume that $\pi$ is minimal in the sense 
that ${\Kc}$ is the smallest invariant subspace for $\pi$ 
containing the range of $V$. 
Now, under these conditions, we can state the 
following theorem which will be very useful for our 
study of Toeplitz algebras.

\proclaim{Theorem~2.2}
Let $\Phi, \Phi _0,\cle,\cstare,{\Kc},V$ and $\pi$ be as above. 
Then $\Ker \Phi _0=\Ker\pi$ and the mapping 
$$\rho\colon \pi(\cstare)\to\bh$$
defined by $\rho(\pi(X))=V^*\pi(X)V$ for 
$X\in\cstare$ is a complete  isometry whose range 
equals $\cle=\Ran \Phi $. 
Moreover, if $\Ran \Phi $ is 
$\sigma$-weakly closed, then $\pi(\cstare)$ is also 
$\sigma$-weakly closed, hence a von Neumann 
subalgebra of $\B({\Kc})$ and the map $\rho$ defined 
above is a $\sigma$-weak homeomorphism. 
\endproclaim

\demo{Proof}
First of all, one can easily see that $\Ker \Phi _0$ is an 
ideal in $\cstare$. 
Indeed, Theorem~2.1 implies that 
$\Ker \Phi _0$ is invariant under multiplication by 
elements in $\cle$ and then use the fact that, since 
$\cle$ is selfadjoint, the \cstar $\cstare$ is the 
closed linear span of all finite products of elements 
from $\cle$ and the identity. 
Now, in order to prove 
the equality of the two kernels, we fix $T\in\cstare$ 
such that $\Phi_0(T)=0$ and let $X,Y\in\cstare$ and 
$\xi,\eta\in {\Hc}$ be arbitrary. 
Then
$$(\pi(T)\pi(X)V\xi,\pi(Y)V\eta)=
(V^*\pi(Y^*TX)V\xi,\eta)=
(\Phi _0(Y^*TX)\xi,\eta)=0$$
because $\Ker \Phi _0$ is an ideal.
Since $\pi$ is minimal, this shows that $\pi(T)=0$. 
Since the other inclusion is trivial, the equality of the 
two  kernels is proved.

We now show that the mapping $\rho$ is a complete 
isometry. First, we see that, since $\Phi\sb 0=\rho\circ\pi$ and 
$Ker\pi=Ker\Phi\sb 0$ it follows that $\rho$ is one-to-one. 
Moreover, since $\Phi _0^2=\Phi _0,$ we have that 
 $\rho\circ\pi\circ\rho\circ\pi=\rho\circ\pi$ hence
 $\pi\circ\rho$ is the identity on $\pi(\cstare).$ It 
then follows, since both $\pi$ and $\rho$ are completely 
contractive, that  $\rho$ is actually completely isometric.    
The last assertion of the theorem follows easily from 
the previous one and the separate weak*-continuity 
of the multiplication on a von Neumann algebra.
 The proof of the theorem is completed.
\qed
\enddemo

As we have mentioned above, this result offers an alternate 
proof of Theorem~3.1 in \cite{CE77}. 

\proclaim{Theorem~2.3}
Let $\Phi \colon \bh\to\bh$ be 
a completely positive and completely contractive   
mapping such that $\Phi ^2=\Phi $. 
Then $\Ran\Phi$ is completely isometric with a unital 
\cstar where the product is defined by the rule 
$$\Phi (X)\circ \Phi (Y)=\Phi (\Phi (X) \Phi (Y))$$
for all $X,Y\in\bh$.
\endproclaim

\demo{Proof}
Using the notation in Theorem~2.2, we see that the 
above defined product is precisely the one induced 
from $\pi(\cstare)$ via the complete isometry 
$\rho$. 
\qed
\enddemo

\remark{Remark~2.4} 
It is a well-known fact that if $A$ is a unital $C^*$-algebra and $\theta:A\to B(\Hc)$ is a unital completely isometric mapping, then there exists a *-homomorphism   
$$\pi:C^*(\theta(A))\to A$$ such that $\pi\circ\theta=id_A,$ 
(see Theorem 4.1 in \cite{CE76}). Therefore, in the case when the mapping $\Phi$ in Theorem 2.2 is unital and  assuming  Theorem 2.3 one can immediately see that the mapping $\Phi_0$ appearing in Theorem 2.2 becomes a *-homomorphism when its range is endowed with the multiplication defined in Theorem 2.3. This offers a shorter proof of Theorem 2.2; however the line we took in that theorem gives simultaneously the isomorphism in Theorem 2.3 and a spatial representation for that algebraic structure, that will be useful in the sequel.    
\qed
\endremark

We shall apply Theorem~2.2 to the study of the $C^*$-
algebra
generated by a commuting family of spherical 
isometries.

\definition{Definition~2.5}
A commuting family $\cals=\pti$ of bounded 
operators on a Hilbert space ${\Hc}$ is said to be a 
{\it spherical isometry} if
$$\sum\sb{j\in J}T^*_jT_j=1$$
in the weak operator topology.
\qed
\enddefinition

For instance, if $m$ is any probability Borel measure on the unit sphere $S^{2n-1}$ in ${\Bbb C}^n,$ and $\Hc\subset L^2(m)$ is a jointly  invariant subspace  for the multiplication operators 
$\{M_{z_1},\dots, M_{z_n}\}$ on $L^2(m),$ then their restrictions 
$\{T_{z_1},\dots, T_{z_n}\}$
to $\Hc$ form a spherical isometry. Of particular interest is the case when $m$ is the normalized area measure and $\Hc$ is the $L^2(m)$-closure of all analytic polynomials (the Hardy space  $H^2(S^{2n-1})$) in which case 
$\{T_{z_1},\dots, T_{z_n}\}$ is called the Szego n-tuple on $H^2(S^{2n-1}).$  

If $\pti$ is an arbitrary family of operators on ${\Hc}$ 
satisfying the above equation in particular a spherical 
isometry   then a completely positive  unital, hence 
completely contractive mapping 
$\phi\colon \bh\to\bh$ can be defined by the formula
$$\phi(X)=\sum\sb{j\in J}T_j^*XT_j$$ 
and is also 
obvious that $\phi$ is  $\sigma$-weakly continuous. 
\par
Our main object of study in what follows is a 
commuting family of spherical isometries 
$\calf=\salf$ on some Hilbert space ${\Hc}$. This means 
that if $\alpha\in\Gamma$ then $\salfa=\{\tja\}\sb{j\in 
J\sb\alpha}$  is a spherical isometry and the union $\cup\sb{\alpha\in\Gamma}\salfa$
 is a commutative  set of operators. 

\definition{Definition~2.6}
Given a commuting family $\calf=\salf$ of spherical 
isometries on ${\Hc}$ we define, using the notations 
above,  the space $\calt(\calf)$ of all 
{\it $\calf$-Toeplitz operators} 
to be the set of all operators 
$X\in\bh$ such that
$$\sum\sb{j\in 
J\sb\alpha}T\sb{j,\alpha}^*XT\sb{j,\alpha}=X$$
for all $\alpha\in\Gamma$.
\qed 
\enddefinition

In other words, $\calt(\calf)$ is the set of all 
common fixed points of the completely positive 
 mappings associated to each 
spherical isometry from $\calf$. It is obvious that 
$\calt(\calf)$ contains the commutant
of $\calf$ in particular it contains all the sets $\salfa$ for 
$\alpha\in\Gamma.$   
We shall construct a completely positive projection 
$\Phi$ on this space which will play a crucial role for 
our study of Toeplitz algebras associated to spherical 
isometries. In order to do this, we need 
the following lemma which is a particular case of a 
more general result proved in \cite{BP05}. 
However for 
completeness we shall give below a direct proof.

\proclaim{Lemma~2.7}
Let $\{\phi\sb\alpha\}\sb{\alpha\in\Gamma}$ be a set 
of commuting completely positive unital and  
$\sigma$-weak continuous mappings acting on 
$\bh$ for some Hilbert space ${\Hc}$.  
Then there exists 
a completely positive mapping $\Phi \colon \bh\to\bh$ whose 
range is precisely the set 
$$\{X\in\bh\colon \phi\sb\alpha(X)=X,\ \alpha\in 
\Gamma \}$$
and such that $\Phi ^2=\Phi $.
\endproclaim

\demo{Proof}
Let $S$ denote the semigroup of all finite products of 
elements from the set 
$\{\phi\sb\alpha\}\sb{\alpha\in \Gamma }$. Each element  
$s\in S$ corresponds to a completely positive unital 
and  $\sigma-$weak continuous mapping 
$\psi\sb s\colon \bh\to\bh$ which is a finite product of 
$\phi\sb\alpha$'s. It is obvious that the fixed point 
set of $\{\phi\sb\alpha\}\sb{\alpha\in \Gamma }$ is the 
same as that of $\{\psi\sb s\}\sb{s\in S}$.
We thus obtain an action
$$\gamma\colon S\times\bh\to\bh$$ defined by 
$$\gamma(s,X)=\psi\sb s(X)$$ for all $s\in S$ and 
$X\in\bh$.
Since $S$ is commutative, a well-known result of 
Dixmier \cite{Di50} shows that $S$ is amenable, which 
means that there exists a state $\mu$ on the $C^*$-algebra
$\ell\sp\infty(S)$ of all bounded complex functions 
on $S$ which is invariant under all translations with 
elements from $S$. More precisely, if $t\in S$ and 
$L\sb t\colon \ell\sp\infty(S)\to\ell\sp\infty(S)$ is 
defined by $L\sb t(f)(s)=f(ts)$ for $s\in S$ then 
$\mu(L\sb tf)=\mu(f)$ for all $f\in\ell\sp\infty(S)$.

Now, given $T\in\bh$, for each pair of vectors 
$\xi,\eta\in {\Hc}$ define 
$[\xi,\eta]\sb T=\mu(\gamma(\cdot,T)\xi,\eta)$
and observe that this is a 
bounded sesquilinear map therefore 
there exists an operator that we shall denote by 
$\Phi (T)$ in $\bh$ such that
$$(\Phi (T)\xi,\eta)=[\xi,\eta]\sb T$$
for all $\xi,\eta\in {\Hc}$.
It is now a matter of routine to verify that the 
mapping $T\mapsto \Phi (T)$ is completely positive. 
It is  straightforward to see that if $T\in\bh$ is such 
that 
$\psi\sb s(T)=T$ for all $s\in S$ then $\Phi(T)=T$ as 
well. 

We will show now that 
$\psi\sb s(\Phi (T))= \Phi (T)$ for all $T\in\bh$. 
In order to see that, recall that all the mappings 
$\psi\sb s$ are $\sigma$-weakly continuous, which 
means that for any $s\in S$ there exists a norm 
continuous mapping 
$\psi\sb s\sp *$ on the space ${\frak S}\sb 1({\Hc})$ 
of trace-class operators on ${\Hc}$ such that
$\tr(\psi\sb s(T)L)=\tr(T\psi\sb s\sp *(L))$ for all 
$T\in \bh$ and $L\in{\frak S}\sb 1({\Hc})$. Moreover one 
can 
see that the mapping $\Phi$ satisfies the identity
$\tr(\Phi (T)L)=\mu(\tr(\gamma(\cdot,T)L)$ for all 
$T\in\bh$ and $L\in{\frak S}\sb 1({\Hc})$.
It then follows that for any $t\in S$ we have
$$\tr(\psi\sb t(\Phi (T))L)=
\tr(\Phi (T)\psi\sb t\sp *(L))=
\mu(\tr(\gamma(\cdot,T)\psi\sb t\sp *(L))=
\mu(\tr(\gamma(t\cdot,T)L))$$
and because $\mu$ is an invariant mean, 
this last term equals
$\mu(\tr(\gamma(\cdot,T)L))$ which is equal to 
$\tr(\Phi (T)L)$
for all $T\in\bh$ and $L\in{\frak S}\sb 1({\Hc})$.  
This shows 
that indeed $\psi\sb t(\Phi (T)= \Phi (T)$ for all $T\in\bh$, 
hence $\Phi ^2=\Phi $ as well. 
The proof of this lemma is completed.
\enddemo

In the proof of our main result we also need the 
following easy lemma.

\proclaim{Lemma~2.8}
Suppose $\calm\subset\bh$ is a von Neumann 
algebra and $\calm'$ denotes its 
commutant. Let $\calr=\{T\sb j\}\sb{j\in J}$ be a 
family of operators in $\calm$ such that for all $X\in\calm$ 
$$\sum\sb{j\in J}T\sb j\sp *XT\sb j=X$$ in the 
$\sigma-$weak topology. 
Then $\calr\subset\calm\cap\calm'.$
\endproclaim

 The main result of this paper is the following.

\proclaim{Theorem~2.9}
Let $\calf=\salf$ be a commuting family of spherical 
isometries on some Hilbert space ${\Hc}$ with 
$\salfa=\{\tja\}\sb{j\in J\sb\alpha}$ for each 
$\alpha\in \Gamma $  and let $\calt(\calf)$ be the space of 
all $\calf$-Toeplitz operators (see {\rm Definition~2.6}
above). 
Let also $\cstartf$ denote the $C^*$-subalgebra of $\bh$ 
generated by $\cltf$. Then we have:
\roster

\item "(1)" There exists a completely positive unital  
mapping $\Phi \colon \bh\to \bh$ such that $\Phi ^2=\Phi $ and 
whose range coincides with $\cltf$. 

\item "(2)" There exist a Hilbert space $\tilh$ 
containing ${\Hc}$  and a commuting family 
$\tilf=\tsalf$ of normal spherical isometries on 
$\tilh$ with $\tsalfa=\{\ttilj\}\sb{j\in 
J\sb\alpha}$ 
which leaves ${\Hc}$ invariant and whose restriction to 
${\Hc}$ coincides with $\calf$ in other words the family 
$\calf$ is subnormal.

\item "(3)" Suppose that the normal extension 
$\tilf$ is minimal, i.e. $\tilh$ is the smallest reducing 
subspace for $\tilf$ containing ${\Hc}$. Then there 
exists a unital *-representation
$$\pi\colon \cstartf\to\btilh$$
such that: 

{\rm (3a)} $\pi(\tja)=\ttilj$ for all $\alpha\in \Gamma $ and $j\in 
J\sb\alpha$.

{\rm (3b)} If $P\sb {\Hc}$ is the orthogonal projection of 
$\tilh$ 
onto ${\Hc}$ then
$$\Phi(X)=P_{\Hc}\pi(X)|\sb {\Hc}$$ for every $ 
X\in\cstartf$.

{\rm (3c)} The image $\pi(\cstartf)$ of $\pi$ coincides with 
the commutant in $\btilh$ of the \cstar 
$C^*(\tilf)$ generated by $\tilf$ in $\btilh$.

{\rm (3d)} The mapping
$$\rho\colon \pi(\cstartf)\to \bh$$
defined by $\rho(\pi(X))=P\sb {\Hc}\pi(X)|\sb {\Hc}$ for 
$X\in \cstartf$ is a complete isometry onto the space 
$\caltf$ of all $\calf$-Toeplitz operators such that
$\pi\circ\rho$ is the identity on $\pi(\cstartf).$  
Therefore the short exact sequence
$$0\to \Ker\pi\hookrightarrow\cstartf 
@>{\pi}>>\pi(\cstartf)\to 0$$
has  a completely isometric cross section.
Moreover, $Ker\pi$  coincides with the closed 
two-sided ideal of $\cstartf$
generated by all operators of the form 
$XY-\Phi(XY)$ with $X,Y\in\caltf.$ 
 
{\rm (3e)}  
If $\cstarf$ denotes 
the unital \cstar   generated by $\calf$ in $\bh$
then  $\Phi(\cstarf)=\cstarf\cap\caltf.$ 
Moreover the  kernel of the restriction of $\pi$ to 
$\cstarf$ coincides with the closed ideal of $\cstarf$ generated by 
all the commutators 
$XY-YX$ with $X,Y\in\caltf\cap\cstarf$
hence it coincides with 
the commutator ideal $\tilc$ of $\cstarf.$
Therefore   we have a short exact sequence
$$0\to\tilc\hookrightarrow\cstarf @>{\pi}>> 
C^*(\tilf)\to 0$$ for which the restriction of 
$\rho$ to $C^*(\tilf)$ is  a completely isometric cross section. 

{\rm (3f)} An operator $X\in\bh$ belongs to the commutant 
of $\calf$ if and only if both $X$ and $X^*X$ belong 
to the space $\caltf$ of $\calf$-Toeplitz operators. 
In this case  there exists  a unique operator $\tilx$ in 
the commutant of $\tilf$ which leaves ${\Hc}$ invariant 
and whose restriction to ${\Hc}$  coincides with $X$. 
Moreover the map $X\mapsto\tilx$ is norm preserving.
\endroster

\endproclaim

\demo{Proof}
For each $\alpha\in \Gamma $ let 
$\phi\sb\alpha\colon \bh\to\bh$ be the completely 
positive  $\sigma$-weakly continuous mapping 
associated to the spherical isometry $\salfa$, so for all 
$X\in\bh$ 
we have 
$\phi\sb\alpha(X)=\sum\limits_{j\in 
J_\alpha}T_{j,\alpha}^* X T_{j,\alpha}$. 
It follows that $\cltf$ is precisely the set of common 
fixed points  of the commuting family of mappings 
$\{\phi\sb\alpha\}\sb{\alpha\in \Gamma }$. 
Therefore we 
can apply Lemma~2.7 to infer the existence of an 
idempotent completely positive  mapping 
$\Phi \colon \bh\to\bh$ whose 
range is precisely $\cltf$. 
This proves item~(1).

Let $\Phi $ as in item~(1) and let $\Phi \sb 0$ denote its 
restriction to $\cstartf$. 
Denote by 
$\pi\colon \cstartf\to \B(\tilh)$  the minimal Stinespring 
dilation of $\Phi \sb 0$. Therefore there exists an 
isometry $V\colon {\Hc}\to\tilh$ such that 
$$\Phi \sb 0(X)=V^*\pi(X)V$$ for all $X\in\cstartf$. 
We  see that we are precisely in the situation of 
Theorem~2.2 above, and moreover the range of $\Phi$ is 
also $\sigma$-weakly closed because it is the set of 
all common fixed points of a family of \sw 
continuous mappings. 
The conclusion  that follows 
from Theorem~2.2 is that the mapping 
$$\rho\colon \pi(\cstartf)\to \bh$$
defined by $\rho(\pi(X))=V^*\pi(X)V$ for $X\in 
\cstartf$ is a complete isometry onto the space of all 
$\calf$-Toeplitz operators and that the image of 
$\pi$ is a von Neumann subalgebra of $\B(\tilh)$. 
Let 
$$\ttilj=\pi(\tja)$$ 
for all $\alpha\in \Gamma $ and 
$j\in J\sb\alpha$ and let also denote 
$\tsalfa=\{\ttilj\}\sb{j\in J\sb\alpha}$ and let 
$\tilf=\tsalf$.
Our next aim  is to show that each family $\tsalfa$ is 
a spherical isometry and that 
$$\smj\ttiljs\pix\ttilj=\pix$$ for all $X\in\cstartf$.
For this purpose, fix $\alpha\in\Gamma$ and observe that 
$\sum\limits_{j\in F}\ttiljs\ttilj\le 1$ for each finite subset 
$F\subset\jalfa$. 
Therefore 
$$\smj\ttiljs\ttilj\le 1.$$ 
Now, let $X\in\cstartf$ and 
let $F\subset\jalfa$ be a finite set. 
Then we see that 
$$\rho(\smf\ttiljs\pix\ttilj)=\smf\tjas\Phi(X)\tja.$$
Taking weak*-limits in both sides, using the fact that 
$\rho$ is a weak*-homeomorphism and using that 
$\rho$ is isometric we infer that 
$$\smj\ttiljs\pix\ttilj=\pix$$ for all 
$X\in\cstartf$. 
In particular it follows that 
$\tsalfa$ is indeed a spherical isometry. 
Moreover, 
using Lemma~2.8 we infer that all $\ttilj$ belong to the 
center of $\pi(\cstartf)$, in particular they are commuting normal 
operators.
Since $\tja=V^*\ttilj V$ and both $\salfa$ and 
$\tsalfa$ are spherical isometries it is easy to see 
that $\ttilj V{\Hc}\subset V{\Hc}$ for all $\alpha\in\Gamma $ and 
$j\in J\sb\alpha$.
This shows that the family $\calf$ is subnormal, 
which proves item~(2).
 
We will show now that $\tilf$ is the minimal normal 
extension of $\calf$. 
For this purpose let ${\Kc}$ be the 
smallest reducing subspace for $\pi(\cstarf)$ 
containing $V{\Hc}$. 
Let $\pik\colon \cstarf\to\bk$ the 
*-representation  defined by $\pik(X)=\pk\pi(X)|_{\Kc}$ 
where $\pk$ denotes the orthogonal projection of 
$\tilh$ onto ${\Kc}$. 
We will show that the map defined by 
$\rhok(\pix)=\pk\pix|_{\Kc}$ is a *-isomorphism of 
$\pi(\cstartf)$ onto the commutant  
$\pik(\cstarf)'$ of $\pik(\cstarf)$. 

It is clear that 
$\rhok$ is a completely positive and completely 
contractive mapping. It takes values in  
$\pik(\cstarf)'$ because each $\ttilj$ is in the center 
of $\pi(\cstartf)$ and because the space ${\Kc}$ is 
reducing for all $\ttilj$. 
Let $\rhoh\colon \pik(\cstarf)'\to\bh$ be defined by 
$\rhoh(Y)=V^*YV$ for $Y\in\pik(\cstarf)'$. 
Then it 
is obvious that its image is included in $\calt(\calf)$ 
and moreover 
$\rho=\rhoh\rhok$ where 
$\rho\colon \pi(\cstartf)\to\bh$ was defined above as 
$\rho(Y)=V^*YV$ for $Y\in\pi(\cstartf)$. 
Recall 
now that we already proved that $\rho$ is 
completely isometric which implies that the mapping 
$\rhok$ is completely isometric. 
Therefore in order 
to show that $\rhok$ is onto, it suffices to show that 
the mapping $\rhoh$ is one-to-one. 
Suppose 
therefore that $Y\in\pik(\cstarf)'$ is such that 
$\rhoh(Y)=V^*YV=0$. 
In order to show that $Y=0$ 
it suffices, because $\pik$ is a minimal dilation of 
$\Phi $ restricted to $\cstarf$ and $\pik(\cstarf)$ is 
abelian and $Y$ is in its commutant, to show that for 
any two finite families $\{\tj1,\dots,\tjm\}$ and 
$\{\ti1,\dots,\tin\}$ we have
$V^*\ttj1^*\dots\ttjm^*Y\tti1\dots\ttin V=0$ and 
the latter equality follows immediately from the fact 
that $V{\Hc}$ is invariant for all $\tilf$.
This shows that $\rhoh$ is indeed one-to-one hence 
$\rhok$ is onto.  
Since, by a well known result 
of Kadison [Kad51], any completely 
isometric surjective unital mapping 
between two $C^*-$ algebras is multiplicative
it follows that 
$\rhok$  is indeed a 
*-isomorphism of $\pi(\cstartf)$ onto  
$\pik(\cstarf)'$ in particular the space ${\Kc}$ is 
invariant under $\pi(\cstartf)$. 
Since $\pi$ is 
minimal, this shows that in fact we have that 
${\Kc}=\tilh$. 
In particular this shows that 
$\pik(X)=\pix$ for all $X\in\cstarf$.   
Moreover, the fact that $Ker\pi$ is 
the ideal generated by all operators of the form
$XY-\Phi(XY)$ with $X,Y\in\caltf$
follows by an easy induction argument 
on the length of an arbitrary product
of elements from $\caltf$ using the fact 
that $Ker\pi=Ker\Phi_0$
which equals $\{X-\Phi(X)\colon X\in\cstartf\}$
together with Theorem 2.1.  
This completes the proof 
of  (3a), (3b), (3c)
and (3d).

In order to prove (3e) we  show first that 
$\Phi(\cstarf)=\cstarf\cap\caltf.$ Since 
$\Phi^2=\Phi$ it is enough to show that 
 $\Phi (\cstarf)\subset\cstarf$. This inclusion 
follows easily from the fact that since 
$\pik(\cstarf)$ is abelian, $\Phi $ takes any finite 
product of $\tja$'s and $\tja^*$'s into a 
permutation of the same  product having all  the 
$\tja^*$'s at the left and all the $\tja$'s  at the 
right.

Now we can easily  prove that the kernel of $\pik$   coincides 
with the  ideal  of $\cstarf$ 
generated by all commutators $XY-YX$ with 
$X,Y\in\cstarf\cap\caltf$.
First, since $\pik(\cstarf)$ is commutative, we have that 
 any such commutator is in $ \Ker(\pik)$. 
 Let us denote by 
$\Phi \sb{00}$ the restriction of $\Phi $ to $\cstarf.$
We see from the proof of Theorem~2.2 that $\Ker \Phi \sb 0
=\Ker\pi$ therefore $\Ker\pik=\Ker \Phi \sb{00}$ as 
well. On the other hand, since $\p00^2=\p00$ we 
see that $\Ker\p00=\Ran(I-\p00)$. 
Now,  if  
$X\in\cstarf$ is  a finite product of  $\tja$'s and $\tja^*$'s
it becomes obvious from the above description of $\p00(X)$ that
$X-\p00(X)$ belongs to the ideal generated by all commutators 
$XY-YX$ with  
$X,Y\in\cstarf\cap\caltf.$ 
This completes the proof of (3e). 

An alternate proof of the fact that $Ker\pik$ coincides with the commutator ideal in $\cstarf$ can be based on Bunce'  characterization of  multiplicative functionals on $C^*$-algebras generated by commuting hyponormal operators in terms of their joint approximate point spectrum (see [Bun71]). Indeed, in our case, it can be easily shown that for any $\{T_1,\dots,T_m\}\subset\calf,$ we have that $\sigma_{ap}(T_1,\dots,T_m)$ equals $\sigma_{ap}(\pi(T_1,\dots,\pi(T_m))$ where $\sigma_{ap}$ stands for the joint approximate spectrum. 

We now  prove (3f).
 If $X\in\bh$ is such that $X$ commutes with all 
operators from ${\Cal F}$ then obviously $X$ and 
$X^*X$ belong to $\calt(\calf)$. 
Suppose now that $X\in\bh$ is such that both $X$ 
and $X^*X$ belong to $\cltf$.
If $\tilx=\pi(X)$ then $\tilx$ commutes with all the 
normal extensions from $\tilf$ and $V^*\tilx V=X$. 
Moreover $\Vert\tilx\Vert=\Vert X\Vert$ so all we 
need to
show is that $\tilx V{\Hc}\subset V{\Hc}$. For this 
purpose, 
we observe that since  $X^*X\in\calt(\calf)$ then 
$$X^*X=V^*\pi(X^*X)V=V^*\tilx^*\tilx V.$$  
Therefore 
if $\xi\in {\Hc}$ then
$$\Vert V^*\tilx V\xi\Vert=
\Vert X\xi\Vert=
\Vert\tilx V\xi\Vert$$
which implies that indeed $\tilx V{\Hc}\subset V{\Hc}$.
This finishes the proof of (3f) and the proof of the 
theorem as well.
\qed
\enddemo

Let us note that both in the  case of the unilateral shift $S$ on $H^2({\Bbb T})$,   and in the more general case of the Szego n-tuple all the assertions from (2) and (3) are well known. We should refer here to the pioneering work of L.A. Coburn on the $C^*$-algebra generated by a single isometry; see  \cite{Co67} and \cite{Co69}.
In particular, in the first case   (3d) and (3e) are classical results due to L. Coburn and R.G. Douglas, and  moreover,  the kernel of $\pi$ at (3d) is the corresponding commutator ideal; see  
Chapter VII in  \cite{Do98} for a self-contained presentation of the classical theory and precise references to the original papers.
In the case of the Szeg\H o  $n$-tuple on the unit sphere in ${\Bbb C}^n$  (3e) is  proved in \cite{Co73}, 
and (3d) appears in \cite{DJ77}. 
We mention also that in both these cases the commutator ideal appearing in (3e) was shown to 
coincide with the ideal of all compact operators on the corresponding 
Hardy space. 

For the case of a commuting family of isometries, 
the existence of a commuting unitary extension was proved 
in \cite{Ito58}. In this case, the commutant lifting at (3f) is quite straightforward once we assume Ito's result. These results hold true even on Banach spaces, see \cite{Do69}. 
Exact sequences similar to that in (3e) for 
finite families of commuting isometries have been studied in \cite{BCL78}.

$C^*$-algebras generated by  isometric representations of 
commuting semigroups  have been studied mainly for semigroups of  positive elements  in ordered abelian groups, see for example \cite{BC70}, \cite{Do72},\cite{Mu87} or more recently  \cite{ALNR94} for a study based on crossed products  by endomorphisms, a line which has been intensively pursued in the last decade.  

For the case of a single finite spherical isometry the existence of 
a normal extension along with 
a commutant lifting theorem  were proved in \cite{At90}; 
see also \cite{AL96} for alternate proofs. 

\head 
3. Toeplitz operators associated with function algebras
\endhead

The following general framework is frequently used when dealing with Toeplitz operators on Hardy spaces. 
Let $K$ be a compact  Hausdorff space and let $\ck$ denote 
the Banach algebra of all complex-valued continuous 
functions on $K.$  Let 
$A\subset\ck$ be a norm-closed subalgebra containing the 
constants and separating the points of $K.$ Such algebras 
are called function algebras or uniform algebras 
(see \cite{Ga69} for basics of uniform algebras). 
Let us  consider a Borel 
probability measure $m$ on $K$ and let $supp(m)$ be its 
closed support. The generalized Hardy space $\har$ 
associated to $A$ is  the $\lar$ closure of $A.$ 
For any function $\phi\in\linf$ the Toeplitz operator 
$\tfi\colon\bhar\to\bhar$ is defined by 
$\tfi h=P_{\har}(\phi h)$ for $h\in\har$ where $P_{\har}$ is the orthogonal 
projection of $\lar$ onto $\har.$ 
We shall also consider the 
usual multiplication operators $\mfi$ defined on $\lar$ by 
$\mfi f=\phi f$ for all $f\in\lar.$ 
Let $\hinfty$ denote the intersection $\har\cap\linf$ which is 
a weak*-closed subalgebra of $\linf.$ 
If  $B\subset\linf$ is any 
unital subalgebra, we shall denote by $\Cal T(B)$ the $C^*-
$subalgebra of $\bhar$ generated by all Toeplitz operators 
$T\sb\phi$ with $\phi\in B$ and by $\Cal C(B)$ 
the closed ideal in $\Cal T(B)$ generated by all operators of 
the form $T\sb\phi T\sb\psi-T\sb{\phi\psi}$ for arbitrary 
$\phi,\psi\in B.$ 

For our purposes we need to introduce the following 
definition. We shall say that a finite family of functions 
$F=\{\phi\sb 1,\dots,\phi\sb n\}\subset\ck$ is a spherical 
multifunction  if 
$$\sum\sb{j=1}\sp n|\phi\sb j(x)|^2=1$$ for  every $x\in 
K.$

We are now able to state the main result of this section.

\proclaim{Theorem 3.1}
Let $K$ be a compact Hausdorff space and
let $A\subset\ck$ be a unital norm-closed  subalgebra. 
Suppose  
 there exists  a family $\falfa$ of spherical multifunctions 
in $\ck$ where each 
$F\sb\alpha$ is of the form 
$F\sb\alpha=\{\phi\sb{j,\alpha}\}\sb{j\in\jalfa}$ with each  
$\phi\sb{j,\alpha}\in A$ and such that for each pair of distinct 
points $x,y\in K$ there exist an index $\alpha\in\Gamma$ 
and an index $j\in J\sb\alpha$ such that 
$\phi\sb{j,\alpha}(x)\neq\phi\sb{j,\alpha}(y).$    Then for 
any   Borel probability measure  $m$ on $K$  the following 
assertions hold true.

{\rm(1)} A bounded operator $X\in\bhar$ is a Toeplitz operator if  
and only if  it satisfies the following  
equations: 
$$\sum\sb{j\in\jalfa}T\sb\fijalfa\sp{*}XT\sb\fijalfa=X$$ for 
all $\alpha\in\Gamma$ (in the  case when $A$ is the disc algebra and $m$ is the Lebesgue measure on the unit circle 
we then retrieve the classical result of Brown-Halmos by 
taking as spherical multifunction $\phi(z)=z$).

{\rm(2)} A bounded operator $X\in\bhar$ is of the form 
$X=T\sb\psi$ for some $\psi\in\hinfty$ if and only if it 
commutes with  $T\sb{\phi\sb{j,\alpha}}$ for all 
$\alpha\in\Gamma$ and all $j\in\jalfa,$ if and only if it commutes with all $T\sb\phi$ with $\phi\in A.$
The map $\psi\mapsto T\sb\psi$ is a Banach algebra isometric 
isomorphism between $\hinfty$ and the
 commutant of the family $\{T\sb\phi\colon\phi\in A\}.$ 
Moreover this commutant  is a maximal 
abelian subalgebra of $\bhar$ and hence, 
for every $\phi\in\hinfty$ we have that $\sigma(T\sb\phi)=\sigma\sb{\hinfty}(\phi).$

{\rm(3)} There exists a short exact sequence of $C^*-$algebras  

$$0\to\calinf\hookrightarrow\taulinf @>{\pi}>>\linf\to 0$$ such that 
$\pi(T\sb\phi)=\phi$ for all $\phi\in\linf.$ 
In particular the spectral inclusion  
$\text{\rm essran}(\phi)\subset\sigma(T\sb\phi)$ holds true (in the case of the unit circle this is the classical theorem of Hartman and Wintner [HW54]) and we also have that 
$$\sigma(T\sb\phi)\subset conv(essran(\phi))$$ where $conv$ denotes the convex hull.

{\rm(4)}  There exists a short exact sequence 
$$0\to\calck\hookrightarrow\tauck @>{\pi}>> C(supp(m))\to 0$$
such that $\pi(T\sb\phi)=\phi$ on $supp(m).$ Moreover, in 
this case $\calck$ coincides always with the closed ideal in 
$\tauck$ generated by all commutators $T\sb\phi T\sb\psi-
T\sb\psi T\sb\phi$ with $\phi,\psi\in\ck.$ 

\endproclaim

\demo{Proof}
Let us denote, for each $\alpha\in\Gamma$ and each 
$j\in\jalfa$ by $\tjalfa$ the Toeplitz operator with symbol 
$\phi\sb{j,\alpha}.$ Since each tuple 
$\{\phi\sb{j,\alpha}\}\sb{j\in\jalfa}$ is a spherical multifunction 
it follows 
easily that in this case 
$\salfa=\{\tjalfa\}\sb{\alpha\in\Gamma}$ is a spherical 
isometry and that $\calf=\{\Cal 
S\sb\alpha\}\sb{\alpha\in\Gamma}$ is a commuting family 
of spherical isometries in $\bhar.$ The separation property  
imposed on these spherical multifunctions implies via the Stone-Weierstrass 
theorem that 
the $C^*-$algebra  generated in $\ck$ by the union of all families 
$F\sb\alpha$ with $\alpha\in\Gamma$ equals $\ck$ itself.
In turn this implies that the set 
$\tilf$ of all the corresponding multiplication operators 
$\mfialfa$ on $\lar$ is the minimal normal extension of 
$\calf.$ Therefore, using Theorem 2.9 we infer that every 
operator $X\in\bhar$ satisfying  
equations  (1) is the compression of a bounded 
operator $Y$ in the commutant of all operators $\mfialfa$  
therefore $Y$ commutes with all multiplication operators 
$M\sb\phi$ with $\phi\in\ck$ which implies that $Y$ itself is 
a multiplication operator with some function $\psi\in\linf$ 
which shows that $X$ is a Toeplitz operator i.e. 
$X=T\sb\psi.$ Conversely, any Toeplitz operator obviously 
satisfies these equations  because $\har$ is invariant for all 
operators $\mfialfa.$ This completes the proof of (1). 
Now, the   proofs of (2), (3) and (4) follow easily from the previous remarks  combined with Theorem 2.9 (the last assertion at (3) follows  from the well-known fact that for any bounded operator $T$ we have that $conv(\sigma(T))\subset W(T)$ with equality for normal operators, where $W(T)$ stands for the  numerical range of $T;$ the case of the unit circle is due to A. Brown and P.R. Halmos \cite{BH63}).
\qed
\enddemo

As a remark, it can be shown, on the same lines of the proof 
of Theorem 2.9, or using the results from [Bun71], that if the \cstar generated by $\hinfty$ 
in $\linf$ coincides with $\linf$ (equivalently if $\hinfty$ 
separates the points in the maximal ideal space of $\linf$) 
then $\calinf$ coincides with the closed ideal of $\taulinf$ 
generated by all commutators $T\sb\phi T\sb\psi-T\sb\psi 
T\sb\phi$ with $\phi,\psi\in\linf$ (for instance this is the 
case when $K$ is the unit circle, $A$ is the disc algebra and 
$m$ is the Lebesgue measure on $K;$ see \cite{Do98}).

We also remark that a  description of the character space of 
the quotient 
$$\taulinf/\calinf$$ 
valid for Hardy spaces over 
any function algebra was given in \cite{Sun87}.

Here follow two general  examples of function algebras 
satisfying the hypotheses of Theorem 3.1.  
We emphasize that this result holds for {\it every} Borel 
probability measure on the corresponding Shilov boundaries. 
Toeplitz operators on such domains have been intensively studied 
in the literature; see \cite{Up84} and \cite{Up96}. 

\example{Example 3.2}
Let $\Omega\subset {\Bbb C}^n$ be a bounded strictly pseudoconvex 
domain and let $A(\Omega)$ be the algebra of all 
continuous functions on its closure  
and holomorphic on $\Omega.$ Let $K=\partial\Omega$ be 
the topological boundary of $\Omega$ which coincides in 
this case with the Shilov boundary of $A(\Omega)$ and let 
$A$ be the set of all restrictions to $\partial\Omega$ of all 
functions from $A(\Omega).$ 
It follows from an embedding theorem for such domains 
(see Theorem 3 in \cite{Lo85})  
that there exist  
a natural number $N>1$ and functions $f\sb1,\dots,f\sb N$ 
in $A(\Omega)$ such that the function 
$F\colon\partial\Omega\to {\Bbb C}^N$ defined by 
$F(x)=(f\sb 1(x),\dots,f\sb N(x))$ is one-to-one and takes 
$\partial\Omega$ into the unit sphere in ${\Bbb C}^N.$
This shows that it is a separating spherical
multifunction for $A$ and hence Theorem 3.1 applies in this case.
In particular this applies to any  bounded domain with $C^2$ boundary in the complex plane. For the case of finitely connected domains in ${\Bbb C}$ with analytic boundary, exact sequences of the form (3) and (4) were constructed in \cite{Ab74}.  
\endexample

\example{Example 4.3}
 Let $\Omega\subset {\Bbb C}^n$ be a bounded symmetric domain containing the origin and such that $e^{i\theta}\zeta\in\Omega$ whenever $\zeta\in\Omega$ and  $\theta\in{\Bbb R}.$ Let 
$A(\Omega)$ be as in the previous example. Let  $K\subset\partial \Omega$ be 
the Shilov boundary of $A(\Omega)$ and let   
$\gamma=max\{\vert\zeta\vert\colon\zeta\in\overline{\Omega}\}.$ It is then known that 
   $K=\{z\in\partial\Omega\colon \vert z\vert=\gamma\}$ 
(see Theorem 6.5 in \cite{Loo77}). 
Therefore the function $F(z)={z/\gamma}$  is an 
imbedding of $K$ into the unit sphere in ${\Bbb C}^n$ hence 
Theorem 3.1 applies in this case as well, taking $A$ to be 
the algebra of all restrictions to $K$ of functions from 
$A(\Omega).$  In particular we obtain that, 
given any Borel probability measure $m$ on $K,$ 
 an operator 
$X\in\bhar$ is a Toeplitz operator if and only if 
$$\sum\sb{j=1}\sp{n}T\sb{z\sb j}^*XT\sb{z\sb 
j}=\gamma^2 X.$$  
\endexample

\head References\endhead

\refstyle{ABCDEF}

\enddocument